\documentclass[a4paper,10pt]{article}

\usepackage{amsthm}

\usepackage{amsmath}

\usepackage{amsfonts}

\usepackage{amssymb}

\usepackage{amsmath,amsthm,amscd,amssymb}

\usepackage{latexsym}

\usepackage{graphicx}

\usepackage[all]{xy}

\theoremstyle{plain}

\newtheorem{thm}{\textbf{Theorem}}

\newtheorem{lem}{\textbf{Lemma}}

\newtheorem{df}{\textbf{Definition}}

\newtheorem{cor}{\textbf{Corollary}}

\newcommand{\A}{\Bbb{A}}

\newcommand{\R}{\Bbb{R}}

\newcommand{\C}{\Bbb{C}}

\newcommand{\Q}{\Bbb{Q}}

\newcommand{\F}{\Bbb{F}}

\newcommand{\Z}{\Bbb{Z}}

\newcommand{\Hom}{\text{Hom}}

\newcommand{\rk}{\text{rk}}

\newcommand{\Aut}{\text{Aut}}

\newcommand{\tr}{\text{tr}}

\newcommand{\im}{\text{im}}

\newcommand{\Gal}{\text{Gal}}

\newcommand{\GL}{\text{GL}}

\newcommand{\SL}{\text{SL}}

\newcommand{\e}{\vspace{1pt}}

\newcommand{\y}{\hspace{6pt}}

\title{{\bf{A proof of the Breuil-Schneider conjecture in the indecomposable case}}}

\author{Claus M. Sorensen}

\begin{document}

\date{}

\maketitle

\begin{abstract}
This paper contains a proof of a conjecture of Breuil and Schneider, on the existence of an invariant norm on 
any locally algebraic representation of $\GL(n)$, with integral central character, whose smooth part is given by a 
generalized Steinberg representation. In fact, we prove the analogue for any connected reductive group $G$. This is done by passing to a 
global setting, using the trace formula for an $\R$-anisotropic model of $G$. The ultimate norm comes from classical $p$-adic modular forms.  
\footnote{{\it{2000 AMS Mathematics Classification}}: 11F33, 11F55, 11F72, 11F80.}
\end{abstract}

%----------------------------------------------------------------------------------------------------------

\section{Introduction}

The $p$-adic Langlands program is still in its infancy. For a $p$-adic field $F$, one anticipates a correspondence between certain 
Galois representations $\rho:\Gal(\bar{\Q}_p/F) \rightarrow \GL_n(\bar{\Q}_p)$ and certain representations $\hat{\pi}$ of $\GL_n(F)$ on 
$p$-adic Banach spaces. See Breuil's survey [Br] from the ICM 2010. This correspondence should somehow be compatible with reduction mod $p$, 
cohomology, and $p$-adic families. This is a (big) theorem for $\GL_2(\Q_p)$, due to the work of many people (Berger, Breuil, Colmez, Paskunas, 
and others). However, beyond this example next to nothing is known. Even $\GL_2(F)$, for fields $F \neq \Q_p$, seems surprisingly hard to deal with. 
Let us return to $\GL_2(\Q_p)$ for a moment, and give more details: We start off with a potentially semistable Galois representation
$$
\rho: \Gal(\bar{\Q}_p/\Q_p)\rightarrow \GL(V)\simeq \GL_2(E),
$$ 
with coefficients in a finite extension $E/\Q_p$. We assume $\rho$ is {\it{regular}}. That is, it has distinct Hodge-Tate
weights $w_1<w_2$. By a standard recipe of Fontaine, to be recalled below, one associates a Weil-Deligne representation $\text{WD}(\rho)$. By the 
classical local Langlands correspondence, its Frobenius-seimisimplification $\text{WD}(\rho)^{F-ss}$ corresponds to an irreducible smooth representation
$\pi'$ of $\GL_2(\Q_p)$ over $E$. We let $\pi=\pi'\otimes |\det|^{-1/2}$ if $\pi'$ is generic (that is, infinite-dimensional). If $\pi'$ is non-generic,
we replace it by $\pi=\pi''\otimes |\det|^{-1/2}$, where $\pi''$ is a certain parabolically induced representation with $\pi'$ as its unique 
irreducible quotient. This is the {\it{generic}} local Langlands correspondence. Note that $\pi$ may be reducible. Now, one attaches to $\rho$ an
admissible unitary Banach space representation $B(\rho)$ of $\GL_2(\Q_p)$ over $E$ satisfying a list of desiderata [Br, p. 8]. Most important for us, 
is that $B(\rho)$ is the completion, relative to a suitable invariant norm, of the locally algebraic representation (at least when $\rho$ is 
irreducible):
$$
B(\rho)^{alg}={\det}^{w_1}\otimes_E \text{Sym}^{w_1-w_2-1}(E^2) \otimes_E \pi.
$$
Moreover, $B(\cdot)$ is compatible with the mod $p$ local Langlands correspondence. 

\medskip

\noindent The Breuil-Schneider conjecture mimics some of this for $\GL_n(F)$. Again, let
$$
\rho: \Gal(\bar{\Q}_p/F) \rightarrow \GL(V)\simeq \GL_n(E)
$$ 
be a potentially semistable Galois representation. With $\rho$, we associate a Weil-Deligne representation $\text{WD}(\rho)$ and a multiset
of integers $\text{HT}(\rho)$ as follows: Pick a finite Galois extension $F'/F$ such that $\rho|_{\Gal(\bar{\Q}_p/F')}$ is semistable. Then 
$$
D=(B_{st}\otimes_{\Q_p} V)^{\Gal(\bar{\Q}_p/F')}
$$
is a free $F_0'\otimes_{\Q_p}E$-module of rank $n$, where $F_0'$ is the maximal unramified subfield of $F'$. The module $D$ comes equipped
with a Frobenius $\phi$, a monodromy operator $N$, such that $N\phi=p\phi N$, and a commuting action of $\Gal(F'/F)$. Moreover, there is an 
admissible filtration of $D_{F'}$ by $\Gal(F'/F)$-invariant $F'\otimes_{\Q_p}E$-submodules, which allows to recover $\rho$. 
Observe that one has a factorization,
$$
\text{$D_{F'}\simeq \prod_{\sigma:F \rightarrow E} D_{F',\sigma}$, $\y$ $D_{F',\sigma}=D_{F'}\otimes_{F'\otimes_{\Q_p}E}(F'\otimes_{F,\sigma}E)$.}
$$
Hence, for each $\sigma$, we are given a filtration $\text{Fil}^i(D_{F',\sigma})$ by $\Gal(F'/F)$-invariant free $F'\otimes_{F,\sigma}E$-submodules.
Admissibility means, intuitively, that the Hodge polygon lies beneath the Newton polygon. More formally, one introduces numbers $t_N(D)$ and
$t_H(D_{F'})$ as in [BS, p. 15]. The former is given purely in terms of $\phi$, the latter in terms of the filtration. One requires that 
$t_H(D_{F'})=t_N(D)$, and that $t_H(D_{F'}')\leq t_N(D')$ for any subobject $D' \subset D$ (with the induced filtration).

\medskip

\noindent {\it{Hodge-Tate numbers}}: For every embedding $\sigma:F \rightarrow E$, the $n$-element multiset $\text{HT}_{\sigma}(\rho)$ contains 
$i \in \Z$ with multiplicity $\rk_{(F'\otimes_{F,\sigma}E)}\text{gr}^i(D_{F',\sigma})$. We label these,
$$
\text{gr}^i(D_{F',\sigma})\neq 0 \Leftrightarrow i \in \text{HT}_{\sigma}(\rho)=\{i_{1,\sigma}\leq \cdots \leq i_{n,\sigma}\}.
$$
We say $\rho$ is regular (at $\sigma$) if all the Hodge-Tate numbers $i_{j,\sigma}$ are distinct. 

\medskip

\noindent {\it{Weil-Deligne representation}}: Forgetting the filtration, the $(\phi,N)$-module $D$ gives rise to $WD(\rho)$ as follows. Choose
an embedding $F_0'\hookrightarrow E$ and consider $D_E=D\otimes_{F_0'\otimes_{\Q_p}E}E$ with the inherited monodromy operator $N$, and $W_F$-action
$$
\text{$r(w)=\phi^{-d(w)}\circ \bar{w}$, $\e$ $w\in W_F$.}
$$
(Here $d(w)$ is the power of arithmetic Frobenius induced by $w$, its image in $\Gal(F'/F)$ is $\bar{w}$, and
$\phi$ is the semilinear Frobenius on $B_{st}$.) Note that $r|_{W_{F'}}$ is unramified. This defines $\text{WD}(\rho)=(r,N,D_E)$, 
a Weil-Deligne representation.

\medskip

\noindent Conversely, suppose we are given a Frobenius-semisimple Weil-Deligne representation $(r,N,D_E)$ of $W_{F}$ over $E$, unramified when restricted 
to $W_{F'}$, and for each $\sigma:F \rightarrow E$ a set of $n$ distinct integers,
$$
i_{1,\sigma}<\cdots<i_{n,\sigma}.
$$ 
When does these data arise from a potentially semistable $\rho$? By [BS, p. 14] we know $(r,N,D_E)$ corresponds to a $(\phi,N)\times \Gal(F'/F)$-module
$D$. What we are asking for, is an admissible filtration $\text{Fil}^i(D_{F',\sigma})$ such that
$$
\text{gr}^i(D_{F',\sigma})\neq 0 \Leftrightarrow i\in \{i_{1,\sigma}< \cdots < i_{n,\sigma}\}.
$$
The Breuil-Schneider conjecture asserts this is the case precisely when some locally algebraic representation $\xi \otimes_E \pi$ (constructed from
the given data) carries an invariant norm. That is, a non-archimedean norm $\|\cdot\|$ such that $\GL_n(F)$ acts unitarily. 

\medskip

\noindent {\it{The algebraic representation $\xi$}}: This is constructed out of the tuples $i_{j,\sigma}$. Let
$$
\text{$a_{j,\sigma}=-i_{n+1-j,\sigma}-(j-1)$, $\y$ $a_{1,\sigma}\leq \cdots\leq a_{n,\sigma}$.}
$$ 
That is, write $i_{j,\sigma}$ in the opposite order, change signs, subtract $(0,1,\ldots,n-1)$. The sequence $a_{j,\sigma}$ is identified with
a dominant weight for $\GL_n$, relative to the lower triangular Borel. We let $\xi_{\sigma}$ be the corresponding irreducible algebraic representation
of $\GL_n$, and $\xi=\otimes \xi_{\sigma}$, viewed as an irreducible algebraic representation of the restriction of scalars $\text{Res}_{F/\Q_p}\GL_n$, 
over $E$.

\medskip

\noindent {\it{The smooth representation $\pi$}}: This is constructed out of $(r,N,D_E)$ via a modified local Langlands correspondence. Let 
$\pi^{\circ}$ be the smooth irreducible representation of $\GL_n(F)$ (over $\bar{\Q}_p$) associated with $(r,N,D_E)$ by the usual unitary 
local Langlands correspondence (after fixing a square root of $q=\#\F_F$), 
$$
(r,N,D_E) \simeq \text{rec}(\pi^{\circ}\otimes|\det|^{(1-n)/2}).
$$
The twist $\pi^{\circ}(\frac{1-n}{2})$ does not depend on the choice of $q^{\frac{1}{2}}$, and can be defined over $E$. 
By the Langlands classification (see [Ku] for a useful survey), $\pi^{\circ}$ is the unique irreducible quotient of a parabolically induced representation,
$$
\text{Ind}_P^G(Q(\Delta_1)\otimes \cdots \otimes Q(\Delta_r))\twoheadrightarrow Q(\Delta_1,\ldots,\Delta_r)\simeq \pi^{\circ}.
$$
Here the induction is normalized. The $Q(\Delta_i)$ are generalized Steinberg representations, built from segments of supercuspidals, $\Delta_i$, ordered
in a suitable way. We define
$$
\pi=\text{Ind}_P^G(Q(\Delta_1)\otimes \cdots \otimes Q(\Delta_r))\otimes|\det|^{(1-n)/2}.
$$
By [BS, p. 16], this $\pi$ can be defined over $E$. Note that $\pi$ may be reducible, and it admits $\pi^{\circ}(\frac{1-n}{2})$ as its unique 
irreducible quotient. Moreover, $\pi \simeq \pi^{\circ}(\frac{1-n}{2})$ exactly when the representation $\pi^{\circ}$ is generic. Also, $\pi$ is always
generic [JS]. This is the so-called generic local Langlands correspondence for $\GL_n$.

\medskip

\noindent We are now in a position to state the conjecture, announced in [BS] and [Br].

\medskip

\noindent {\bf{The Breuil-Schneider conjecture}}. {\it{Fix data $(r,N,D_E)$ and $i_{j,\sigma}$ as above, and let $\pi$ and $\xi$ be
the representations constructed therefrom. Then the following two conditions are equivalent,

\begin{enumerate}
\item[(1)] The data arises from a potentially semistable Galois representation. 
\item[(2)] The representation $\xi\otimes_E \pi$ admits a $\GL_n(F)$-invariant norm $\|\cdot\|$.
\end{enumerate}
}}

\medskip 

\noindent The implication $(2) \Rightarrow (1)$ is completely known. A few cases were worked out in [BS], and Hu proved it in general in [Hu]. In fact, 
Hu proves a lot more. He shows that $(1)$ is {\it{equivalent}} to what he refers to as the Emerton condition, which is a purely group theoretic 
statement: With $V$ denoting the space $\xi\otimes_E \pi$, 
$$
\text{\it{(3)} $V^{N_0,Z_M^+=\chi}\neq 0 \Rightarrow |\delta_P^{-1}(z)\chi(z)|\leq 1$,}
$$
for all $z \in Z_M^+$. The implication $(2) \Rightarrow (3)$ is an easy exercise.  

\medskip

\noindent We are concerned with the converse, $(1) \Rightarrow (2)$. Our main result is:
 
\medskip

\noindent {\bf{Theorem A}}. {\it{The conjecture holds when $(r,N,D_E)$ is indecomposable.}}

\medskip

\noindent Recall that indecomposable Weil-Deligne representations are precisely those obtained as follows: Starting with an irreducible 
representation $\tilde{r}:W_F \rightarrow \GL(\tilde{D})$, with open kernel, and a positive integer $ \in \Z_{>0}$, let
$$
\text{$D=\tilde{D}^{\oplus s}$, $\y$ $r=\tilde{r}\oplus \tilde{r}(1)\oplus \cdots\oplus \tilde{r}(s-1)$, $\y$ 
$N:\tilde{r}(i-1)\overset{\sim}{\rightarrow} \tilde{r}(i)$.}
$$
Here $\tilde{r}(i)$ denotes twisting $\tilde{r}$ by the $i$th power of $|\cdot|$, the absolute value on $W_F$, transferred from $F^*$ via 
the reciprocity map. Under the (classical) local Langlands correspondence, $\tilde{D}$ corresponds to a supercuspidal $\tau$, and $D$ corresponds
to the generalized Steinberg representation $Q(\Delta)$, where $\Delta$ is the segment
$$
\Delta=\tau \otimes \tau(1) \otimes \cdots\otimes \tau(s-1).
$$
The Jacquet modules of $Q(\Delta)$ can be made explicit, see Lemma 3.1 in [Hu], for example. They are irreducible if nonzero. From that, it is 
easy to see that condition (3) just amounts to saying $\xi\otimes_E \pi$ has integral central character. In fact, this
was already observed in Proposition 5.3 in [BS], where they also state the resulting conjecture explicitly (as Conjecture 5.5), which is what we prove. 
Our methods work for any connected reductive group $G$ defined over $\Q_p$.  

\medskip

\noindent {\bf{Theorem B}}. {\it{Let $G$ be a connected reductive group over $\Q_p$. Let $\xi$ be any irreducible algebraic representation of 
$G_{\bar{\Q}_p}$, and $\pi$ be any essentially discrete series representation of $G$. Then $\xi\otimes \pi$ admits a $G$-invariant 
norm if and only if its central character
is integral.}}

\medskip

\noindent Taking $G=\text{Res}_{F/\Q_p}\GL(n)$, yields Conjecture 5.5 in [BS]. Indeed, the generalized Steinberg representations coincide with
the essentially discrete series representations, for $\GL(n)$. This Theorem, and its proof, is purely group-theoretical. There is no mention of Galois
representations, and much of the previous discussion is meant to be motivation only.

\medskip

\noindent The proof of Theorem B (which implies Theorem A) is by passing to a global setting, and making use of algebraic modular forms. 
By some sort of averaging over finite (cohomology) groups, we first reduce to the case where $G$ is simple and simply connected, in which case the 
condition on the central character is vacuous. For such $G$, a result of Borel and Harder allows us to find a global model $G_{/\Q}$ such that
$G(\R)$ is compact. If $\pi$ is a discrete series, a trace formula argument (due to Clozel in greater generality) shows that $\xi\otimes \pi$ admits
an automorphic extension. Fixing an isomorphism $\iota:\C \rightarrow \bar{\Q}_p$, we infer that $\pi^K$ sits as a submodule of $\mathcal{A}_{G,\xi}^K$,
a space of classical $p$-adic modular forms. Therefore, $\xi\otimes \pi$ contributes to the direct limit of all
$\xi\otimes\mathcal{A}_{G,\xi}^K$, which in turn embeds in $\mathcal{C}_{G}$, the space of all continuous functions 
$G(\Q)\backslash G(\A_f)\rightarrow\bar{\Q}_p$. This latter space carries a supremum-norm, which is obviously invariant under $G(\A_f)$.

\medskip

\noindent I would like to thank Michael Harris and Marko Tadic, for supplying some psychologically comforting facts about Galois conjugates of
discrete series representations. Moreover, thanks are due to Christophe Breuil, Dinakar Ramakrishnan, and Chris Skinner for their encouragement.

\section{Modular forms on definite reductive groups}

\subsection{The complex case}

\subsubsection{Notation} 

For now, we will study automorphic forms on an arbitrary connected reductive group $G$ over $\Q$ such that $G^{\text{der}}(\R)$ is compact.
Here $G^{\text{der}}$ is the derived subgroup, which is then necessarily an $\R$-anisotropic semisimple group. As is standard, 
$A_G$ denotes the maximal $\Q$-split central torus in $G$, and we choose any central torus $Z_G$ (over $\Q$) containing $A_G$. We will often 
take it to be the whole identity component of the center. $K_{\infty}$ is the maximal compact subgroup of $G(\R)$, which is unique, and possibly 
bigger than $G^{\text{der}}(\R)$. 

\subsubsection{Classical automorphic forms}  

Let $\A=\R\times \A_f$ be the ring of rational adeles. Inside $G(\A)$, we introduce the normal subgroup $G(\A)^1$ cut out by all $|\chi|$, where
$\chi$ ranges over the $\Q$-characters of $G$. It contains $G(\Q)$ as a cocompact discrete subgroup, and one has a decomposition 
$$
G(\A)=A_G(\R)^+\times G(\A)^1.
$$
Automorphic forms are affiliated with a central character, which we fix throughout. That is, we pick an arbitrary continuous 
(possibly non-unitary) character
$$
\omega: Z_G(\Q)\backslash Z_G(\A)\rightarrow \C^*,
$$
and consider the Hilbert space $L_G^2(\omega)$ of all measurable $\omega$-central functions 
$$
\text{$f:G(\Q)\backslash G(\A)\rightarrow \C$, $\y$ $\int_{G(\Q)\backslash G(\A)^1}|f(x)|^2dx<\infty$.}
$$
The right regular representation of $G(\A)$ is completely reducible, and $L_G^2(\omega)$ breaks up into (irreducible) automorphic representations
$\pi=\pi_{\infty}\otimes \pi_f$, each occurring with finite multiplicity $m_G(\pi)$. The space of automorphic forms $\mathcal{A}_G(\omega)$, is the
dense subspace of smooth functions $f$ satisfying the usual finiteness properties under the action of $K_{\infty}$, and the center of the universal
enveloping algebra at infinity. We will restrict ourselves to {\it{algebraic}} $\pi$. That is, we will assume $\pi_{\infty}$ is the restriction of
an irreducible algebraic (finite-dimensional) representation
$$
\xi: G_{\C}\rightarrow \GL(W),
$$
which we fix throughout. Its isotypic component is $\xi\otimes \mathcal{A}_G(\omega)$, where we let

\begin{df}
$\mathcal{A}_{G,\xi}(\omega)=\Hom_{G(\R)}(\xi,\mathcal{A}_G(\omega))=(\xi^{\vee}\otimes\mathcal{A}_G(\omega))^{G(\R)}$.
\end{df}

\noindent This is an admissible smooth representation of $G(\A_f)$, which breaks up as a direct sum $\oplus_{\pi} m_G(\pi)\pi_f$, summing over 
automorphic $\pi$, of central character $\omega_{\pi}=\omega$, such that $\pi_{\infty}=\xi$. We view elements of $\mathcal{A}_{G,\xi}(\omega)$ as
vector-valued functions.

\begin{lem}
As a $G(\A_f)$-module, $\mathcal{A}_{G,\xi}(\omega)$ can be identified with the space of all $\omega_f$-central smooth functions 
$$
\text{$f:G(\A_f)\rightarrow W^{\vee}$, $\y$ $f(\gamma_fx)=\xi^{\vee}(\gamma_{\infty})f(x)$, $\y$ $\forall \gamma \in G(\Q)$.}
$$
\end{lem}

\noindent {\it{Proof}}. One introduces a third space, consisting of all smooth $\omega$-central functions
$$
\text{$f:G(\Q)\backslash G(\A)\rightarrow W^{\vee}$, $\y$ $f(xg)=\xi^{\vee}(g)^{-1}f(x)$, $\y$ $\forall g \in G(\R)$.}
$$
Such a function $f$ gives a $G(\R)$-map $\xi\rightarrow\mathcal{A}_G(\omega)$ by sending a vector $w\in W$ to the automorphic form 
$g \mapsto \langle f(g),w\rangle$. On the other hand, restriction to $G(\A_f)$ identifies it with the space of functions in the lemma. $\square$

\medskip

\noindent {\it{Remark}}. We always assume $\xi$ and $\omega$ are compatible, that is $\omega_{\infty}=\xi|_{Z_G(\R)}$.

\medskip

\noindent By smoothness, as $K$ varies over all compact open subgroups of $G(\A_f)$, one has 
$$
\mathcal{A}_{G,\xi}(\omega)=\underset{K}{\varinjlim} \mathcal{A}_{G,\xi}(\omega)^K,
$$
where $\mathcal{A}_{G,\xi}(\omega)^K$ is the subspace of $K$-invariants, a module for the Hecke algebra $\mathcal{H}_{G,K}$ of all
$K$-biinvariant compactly supported $\C$-valued functions on $G(\A_f)$. Again, for this subspace to be nonzero, we need $K$ and $\omega$ to be 
compatible, in the sense that $\omega_f$ is trivial on $Z_G(\A_f)\cap K$.

\medskip

\noindent {\it{Example}}. When $\xi=1$, we are just looking at the space $\mathcal{A}_{G,1}(\omega)$ of all $\omega_f$-central smooth 
$\C$-valued functions on the profinite (hence compact) set
$$
\text{$\tilde{S}=G(\Q)\backslash G(\A_f)=\underset{K}{\varprojlim} S_K$, $\y$ $S_K=G(\Q)\backslash G(\A_f)/K$.}
$$ 
Moreover, $\mathcal{A}_{G,1}(\omega)^K$ is the space of $\omega_f$-central functions on the finite set $S_K$.

\subsection{The $p$-adic case}

\subsubsection{Notation}

We fix a prime number $p$, an algebraic closure $\bar{\Q}_p$, together with an (algebraic) isomorphism $\iota:\C \overset{\sim}{\rightarrow}
\bar{\Q}_p$. We will occasionally make use of an algebraic closure $\bar{\Q}$, always assumed to be endowed with an embedding
$\iota_{\infty}: \bar{\Q} \hookrightarrow \C$. Correspondingly, $\iota_p=\iota\circ \iota_{\infty}$ is an embedding 
$\bar{\Q}\hookrightarrow \bar{\Q}_p$. Via $\iota$, we base change $\xi$ to an algebraic representation over $\bar{\Q}_p$,
$$
\text{$\iota\xi: G_{\bar{\Q}_p}\rightarrow \GL(\iota W)$, $\y$ $\iota W=W \otimes_{\C,\iota}\bar{\Q}_p$.}
$$
Our central character $\omega$ has a $p$-adic avatar, the continuous character
$$
\text{$\omega_{f,p}: Z_G(\Q)\backslash Z_G(\A_f)\rightarrow \bar{\Q}_p^*$, $\y$ $\omega_{f,p}(z)=\iota\omega_{\xi}(z_p)\cdot\iota\omega_f(z)$.}
$$

\subsubsection{Classical $p$-adic automorphic forms}

All constructions of the previous section can be transferred to $\bar{\Q}_p$ via $\iota$. When we put an $\iota$ in front, we mean
tensoring by $\bar{\Q}_p$, as in $\iota W=W \otimes_{\C,\iota}\bar{\Q}_p$.

\begin{lem}
As a $G(\A_f)$-module, $\iota\mathcal{A}_{G,\xi}(\omega)^K$ can be identified with the space of all $\omega_{f,p}$-central functions 
(smooth away from $p$)
$$
\text{$f:G(\Q)\backslash G(\A_f)\rightarrow \iota W^{\vee}$, $\y$ $f(xk)=\iota\xi^{\vee}(k_p)^{-1}f(x)$, $\y$ $\forall k \in K$.}
$$
\end{lem}

\noindent {\it{Proof}}. Given a complex form $f$, as in the previous lemma, one associates the function
$x \mapsto \iota\xi^{\vee}(x_p)^{-1}\iota f(x)$. It is easy to check that one can recover $f$. $\square$

%\medskip

\begin{df}
$\mathcal{C}_{G}(\omega)=\{\text{continuous $\omega_{f,p}$-central $G(\Q)\backslash G(\A_f)\overset{f}{\rightarrow} \bar{\Q}_p$}\}$.
\end{df}

\noindent Any function $f$, as in the lemma, yields a $K$-map $\iota\xi\rightarrow \mathcal{C}_{G}(\omega)$ by sending $w \in \iota W$ to the 
continuous (in fact, locally algebraic) function $g \mapsto \langle f(g),w\rangle$, and vice versa. Here $K$ acts on $\iota\xi$ through 
the projection to $G(\Q_p)$. We have shown, 
$$
\iota\mathcal{A}_{G,\xi}(\omega)^K=\Hom_K(\iota\xi,\mathcal{C}_{G}(\omega))=(\iota\xi^{\vee}\otimes\mathcal{C}_{G}(\omega))^K.
$$
Note that the image of $K$ in $G(\Q_p)$ is compact open, hence Zariski dense, so that $\iota\xi$ is an irreducible representation of $K$. Let us look 
at its isotypic subspace $\mathcal{C}_{G}(\omega)[\iota\xi]$. That is, the sum of all $K$-stable subspaces isomorphic to $\iota\xi$. This is a 
semisimple $K$-representation, and $\Hom_K(\iota\xi,\mathcal{C}_{G}(\omega))$ is its multiplicity space,
$$
\iota\xi \otimes \iota\mathcal{A}_{G,\xi}(\omega)^K \overset{\sim}{\longrightarrow} \mathcal{C}_{G}(\omega)[\iota\xi]\subset \mathcal{C}_{G}(\omega).
$$
As $K$ varies, these identifications are compatible with inclusions among the spaces $\mathcal{A}_{G,\xi}(\omega)^K$. Taking the direct limit, we end 
up with the injection
$$
\underset{K}{\varinjlim} \iota\xi \otimes \iota\mathcal{A}_{G,\xi}(\omega)^K \hookrightarrow \mathcal{C}_{G}(\omega).
$$
It can be checked that this map is $G(\A_f)$-equivariant. The image is the subspace of locally $\xi$-algebraic functions.
Altogether, we arrive at our key result:

\begin{thm}
There is an injective $G(\A_f)$-map $\iota\xi\otimes\iota\mathcal{A}_{G,\xi}(\omega)\hookrightarrow\mathcal{C}_{G}(\omega)$.
\end{thm}

\subsubsection{Existence of invariant norms}

The space $\mathcal{C}_{G}(\omega)$, being a subspace of $\mathcal{C}(\tilde{S},\bar{\Q}_p)$, has a natural sup-norm,
$$
\|f\|={\sup}_{x\in G(\A_f)}|f(x)|_p={\max}_{x\in G(\A_f)}|f(x)|_p,
$$
which is obviously invariant under the $G(\A_f)$-action, that is $\|g\cdot f\|=\|f\|$.

\begin{cor}
If $\pi=\xi\otimes \pi_f$ is an automorphic representation of $G(\A)$, then
$\iota\xi\otimes \iota\pi_f$ has a natural $G(\A_f)$-invariant norm. (Here $G(\A_f^p)$ acts through $\iota\pi_f^p$, and $G(\Q_p)$ acts
diagonally.) 
\end{cor}

\noindent Since $\iota\xi\otimes \iota\pi_f=(\iota\xi\otimes\iota\pi_p)\otimes \iota\pi_f^p$, we deduce:

\begin{cor}
If $\pi_p$ is an irreducible admissible representation of $G(\Q_p)$, which extends to an automorphic representation of $G(\A)$ of weight $\xi$, then
$\iota\xi\otimes \iota\pi_p$ has a $G(\Q_p)$-invariant norm.
\end{cor}

\noindent This norm is far from canonical. There may be many ways to extend $\pi_p$.

\section{A Grunwald-Wang type theorem}

\subsection{The Grunwald-Wang theorem for $\GL(1)$} 

We briefly recall, from [AT, p. 103], the following result of Grunwald (as corrected by Wang).

\begin{thm}
Given a number field $F$, a finite set of places $S$, and for each $v \in S$ a character $\chi_v$ of $F_v^*$ of finite order, 
there exists a finite order Hecke character $\chi$ of $F$ extending $\chi_S=\otimes_{v\in S}\chi_v$. 
\end{thm}

\noindent Furthermore, the order of $\chi$ can be taken to be the least common multiple of the orders of the $\chi_v$, unless a special case
occurs (where the order of $\chi$ becomes twice that). Given an arbitrary $\chi_S$, we see that it can be extended to a Hecke character 
conditionally: Precisely when some twist $\chi_S|\cdot|_S^{s}$ is of finite order. This is a constraint among the $\{\chi_v\}_{v\in S}$
(as $s\in \C$ depends only on $S$).

\subsection{Clozel's argument on limit multiplicities}

We will use the trace formula in its absolute simplest form. Namely, we will assume, for a moment, that $G$ is {\it{semisimple}}. We 
keep all other assumptions. In particular, $G(\R)$ is compact. The trace formula for $G$ is the following identity,
$$
\tr(\phi:L_G^2)=\sum_{\pi}m_G(\pi)\tr\pi(\phi)=\sum_{\{\gamma\}}\text{vol}(G_{\gamma}(\Q)\backslash G_{\gamma}(\A))O_{\gamma}(\phi),
$$
valid for any test function $\phi \in \mathcal{C}_c^{\infty}(G(\A))$. On the spectral side, we are summing over all automorphic
representations $\pi$. On the geometric side, the sum ranges over $\gamma \in G(\Q)$, up to conjugacy. We denote by $G_{\gamma}$ its stabilizer, and
by $O_{\gamma}$ the orbital integral. Measures are chosen compatibly. 

\medskip

\noindent We wish to quickly outline an argument of Clozel, giving an analogue of the Grunwald-Wang theorem for $G$. We start off with a finite 
set of places $S$ of $\Q$, which we assume contains $\infty$. At each $v\in S$, we are given a discrete series representation $\pi_v^{\circ}$ of
$G(\Q_v)$ (that is, its matrix coefficients are square-integrable). 

\begin{thm}
There is a function $\phi_v^{\circ} \in \mathcal{C}_c^{\infty}(G(\Q_v))$ such that, for every tempered irreducible admisisble representation
$\pi_v$, 
$$
\tr\pi_v(\phi_v^{\circ})=\begin{cases}1, \y \pi_v=\pi_v^{\circ}\\ 0, \y \pi_v\neq \pi_v^{\circ}\end{cases}
$$
(Such a $\phi_v^{\circ}$ is called a pseudo-coefficient of $\pi_v^{\circ}$.)
\end{thm}

\noindent {\it{Proof}}. For $v=\infty$ this is in [CD]. The case $v \neq \infty$ is in [C, p. 278]. $\square$

\medskip

\noindent {\it{Note}}. There may be non-tempered $\pi_v$, for which $\tr\pi_v(\phi_v^{\circ})\neq 0$, but only finitely many. See [C, p. 269] and
[C, p. 280]. Let us introduce $\phi_S^{\circ}=\otimes_{v\in S}\phi_v^{\circ}$. Then $\tr\pi_S(\phi_S^{\circ})\neq 0$ for only finitely many 
representations $\pi_S^{\circ}=\pi_{S,0},\ldots,\pi_{S,r}$. 

\medskip

\noindent With this choice of $\phi_S^{\circ}$, the spectral side becomes
$$
\sum_{\pi^S}m_G(\pi_S^{\circ}\otimes \pi^S)\tr\pi^S(\phi^S)+
\sum_{i=1}^r\sum_{\pi^S}m_G(\pi_{S,i}\otimes \pi^S)\tr\pi_{S,i}(\phi_S^{\circ})\tr\pi^S(\phi^S)
$$
for all $\phi^S \in \mathcal{C}_c^{\infty}(G(\A^S))$. We will take this $\phi^S$ to be of the following form:
$$
\phi^S=\text{vol}(K^S)^{-1}\cdot\text{char}_{K^S},
$$
where $K^S \subset G(\A^S)$ is a compact open subgroup, to be varied. With this choice, the spectral side turns into
$$
\dim\Hom_{G(\Q_S)}(\pi_S^{\circ},(L_G^2)^{K^S})+\sum_{i=1}^r \dim\Hom_{G(\Q_S)}(\pi_{S,i},(L_G^2)^{K^S})\tr\pi_{S,i}(\phi_S^{\circ})
$$
In some sense, the key ingredient of Clozel's proof is the following limit multiplicity formula, based on a method of DeGeorge-Wallach.

\begin{lem}
$\lim_{K^S\rightarrow 1}\text{vol}(K^S)\dim\Hom_{G(\Q_S)}(\pi_{S,i},(L_G^2)^{K^S})=0$ for $i>0$.
\end{lem}

\noindent {\it{Proof}}. This is (a weak version of) Lemma 8, [C, p. 274]. $\square$

\medskip

\noindent Now, let us focus on the geometric side,
$$
\sum_{\{\gamma\}}\text{vol}(G_{\gamma}(\Q)\backslash G_{\gamma}(\A))O_{\gamma_S}(\phi_S^{\circ})O_{\gamma^S}(\phi^S).
$$
Here, by Lemma 5 in [C, p. 271], for sufficiently small $K^S$, the factor $O_{\gamma^S}(\phi^S)=0$ unless $\gamma$ is unipotent. 
Since $G$ is $\Q$-anisotropic, this means $\gamma=1$. In the limit, as
$K^S \rightarrow 1$, the geometric side reduces to just one term,
$$
\text{vol}(G(\Q)\backslash G(\A))\phi_S^{\circ}(1)\text{vol}(K^S)^{-1}.
$$
Here $\phi_S^{\circ}(1)=d(\pi_S^{\circ})>0$ is the formal degree, by the Plancherel formula. See Lemma 9 and 12 in [C]. Putting all this
together, we arrive at the following limit formula,

\begin{thm}
$\text{vol}(K^S)\dim\Hom_{G(\Q_S)}(\pi_S^{\circ},(L_G^2)^{K^S})\underset{K^S\rightarrow 1}{\longrightarrow} 
\text{vol}(G(\Q)\backslash G(\A))d(\pi_S^{\circ})$.
\end{thm}

\noindent This is a weak version of Theorems 1A and 1B in [C], which control ramification away from just one prime. We will not need this. 
On the other hand, Clozel's theorems give lower bounds for $\liminf_{K^S\rightarrow 1}$, not exact limits.

\medskip

\noindent What will be crucial for the applications we have in mind later on, is the following extension theorem, in the vein of
Grunwald-Wang,

\begin{cor}
Let $G$ be a semisimple anisotropic $\Q$-group. Given a discrete series
representation $\pi_S^{\circ}$ of $G(\Q_S)$, where $S$ is a finite set of places of $\Q$, there is an automorphic 
representation $\pi$ of $G(\A)$ such that $\pi_S=\pi_S^{\circ}$.
\end{cor}

\section{Invariant norms on discrete series}

\subsection{Forms of algebraic groups}

We will quote (and use) a result of Borel and Harder on locally prescribed forms of algebraic groups. Recall, if $G$ is an algebraic group over a 
field $F$, an $F$-form of $G$ is an $F$-group $G'$ isomorphic to $G$ over the algebraic closure $\bar{F}$. This gives rise to
a cocycle $c:\Gal(\bar{F}/F)\rightarrow \Aut(G)$ in the obvious way, and identifies the set of equivalence classes of forms with the non-abelian 
Galois cohomology set,
$$
H^1(F,\Aut(G)).
$$
We will take $F$ to be a number field. For each place $v$ of $F$, there is an obvious restriction map
$$
H^1(F,\Aut(G))\rightarrow H^1(F_v,\Aut(G)), 
$$
which on forms corresponds to extending scalars $G' \rightsquigarrow G_v'=G'\otimes_F F_v$.  

\begin{thm}
Let $F$ be a number field, $S$ a finite set of places of $F$, and $G$ an (absolutely) almost simple $F$-group which is either simply connected
or of adjoint type. Then the canonical restriction map is surjective,
$$
H^1(F,\Aut(G))\twoheadrightarrow \prod_{v\in S}H^1(F_v,\Aut(G)).
$$
In other words, given an $F_v$-form $G_v'$ for each $v \in S$, there is an $F$-form $G'$ equivalent to $G_v'$ at places in $S$.
\end{thm}

\medskip

\noindent {\it{Proof}}. This is Theorem B in [BH]. $\square$

\medskip

\noindent If $v$ is a real (infinite) place of $F$, there is always a unique compact form $G_v'$, up to equivalence. The 
corresponding cocycle $c$ is essentially given by the Cartan involution. We immediately deduce the following existence result, which will be 
used in the next section.     

\begin{cor}
Let $G$ be an almost simple $\Q_p$-group which is either simply connected
or of adjoint type. Then there is a model over $\Q$, still denoted by $G$, such that $G(\R)$ is compact.
\end{cor}

\medskip

\noindent {\it{Proof}}. The group $G_{\bar{\Q}_p}\simeq G_{\C}$ has a split model over $\Q$ (even over $\Z$, this is the theory of Chevalley groups),
which we will denote by $G^*$. We apply the Theorem to this group, with $S=\{\infty,p\}$. At $\infty$ we take the compact form of $G_{\R}^*$, at
$p$ we take $G$. $\square$

\subsection{The simple case}

The following result is at the heart of our method.

\begin{lem}
Let $G$ be an almost simple $\Q_p$-group which is either simply connected
or of adjoint type. Let $\xi$ be any irreducible algebraic representation of $G_{\bar{\Q}_p}$, and
$\pi$ be any discrete series representation of $G(\Q_p)$ (both over $\bar{\Q}_p$). Then the locally 
algebraic representation $\xi\otimes \pi$ carries a norm, which is invariant under the $G(\Q_p)$-action.  
\end{lem}

\noindent {\it{Proof}}. The key is to embed this in a global situation. Thus, as in the previous Corollary, we first find
a $\Q$-model $G$ such that $G(\R)$ is compact. With a choice of an isomorphism
$\iota:\C \rightarrow \bar{\Q}_p$, we can confuse $\xi$ and $\pi$ with representations over $\C$ (of $G_{\C}$ and $G(\Q_p)$ respectively). 
We will change notation, and denote the previous $\pi$ by $\pi_p^{\circ}$. Also, we let $\pi_{\infty}^{\circ}=\xi|_{G(\R)}$. Both are discrete series,
so by Corollary 3 there is an automorphic representation $\pi$ of $G(\A)$ such that $\pi_{\infty}=\xi$ and $\pi_p=\pi_p^{\circ}$. By Corollary 2,
we see that $\iota\xi\otimes \iota\pi_p^{\circ}$ has an invariant norm. $\square$ 

\subsection{The semisimple case}

From the simple case, we derive the semisimple case,

\begin{lem}
Let $G$ be a connected semisimple $\Q_p$-group. Let $\xi$ be any irreducible algebraic representation of $G_{\bar{\Q}_p}$, and
$\pi$ be any discrete series representation of $G(\Q_p)$ (both over $\bar{\Q}_p$). Then the locally 
algebraic representation $\xi\otimes \pi$ carries a norm, which is invariant under the $G(\Q_p)$-action.  
\end{lem}

\noindent {\it{Proof}}. Now, suppose $G$ is any connected semisimple $\Q_p$-group, and let $G^{\text{sc}}\twoheadrightarrow G$ be its universal covering
over $\Q_p$, see [PR]. The kernel $\pi_1(G)$ is finite. Being simply connected, $G^{\text{sc}}$ is an actual direct product 
$G_1\times\cdots \times G_r$, of finitely many simply connected simple groups $G_i$. By the main theorem of [Si], the restriction of $\pi$ to 
$G^{\text{sc}}$ is a direct sum of finitely many irreducible admissible representations,
$$
\pi|_{G^{\text{sc}}}\simeq \oplus_{j=1}^s (\tau_{1,j}\otimes \cdots\otimes\tau_{r,j}),
$$ 
where $\tau_{i,j}$ is a discrete series representation of $G_i(\Q_p)$. The restriction $\xi|_{G^{\text{sc}}}$ remains irreducible,
and we continue to denote it simply by $\xi$. It factors as a tensor product $\xi_1\otimes \cdots \otimes \xi_r$, where $\xi_i$ is an irreducible
algebraic representation of $G_{i,\bar{\Q}_p}$. According to Lemma 4, each
$\xi_i\otimes \tau_{i,j}$ has a norm $\|\cdot\|_{i,j}$, invariant under the action of $G_i(\Q_p)$. On the tensor product, where $j$ is fixed for now,
$$
(\xi_1\otimes \tau_{1,j})\otimes \cdots\otimes (\xi_r\otimes \tau_{r,j}),
$$
we put the tensor product norm, see [Sc, p. 110] and Proposition 17.4 therein. 
It has the property that
$$
\|v_1\otimes \cdots\otimes v_r\|_j=\|v_1\|_{1,j}\cdots\|v_r\|_{r,j},
$$
with $v_i \in \xi_i\otimes \tau_{i,j}$. It is defined, for sums of pure tensors, by the formula 
$$
\|v\|_j=\inf\{\max\|v_1\|_{1,j}\cdots\|v_r\|_{r,j} :v=\sum v_1\otimes \cdots\otimes v_r  \}.
$$
Here the maximum is over the same index set as the summation. The infimum is over all possible expressions for $v$.
This tensor product norm $\|\cdot\|_j$ is clearly invariant under $G^{\text{sc}}(\Q_p)$. Taking the maximum of all these, over $j=1,\ldots,s$,
we have constructed a $G^{\text{sc}}(\Q_p)$-invariant norm $\|\cdot\|$ on $\xi\otimes \pi$. Now, to make it invariant under $G(\Q_p)$, we note that
$$
G(\Q_p)/\im(G^{\text{sc}}(\Q_p)\rightarrow G(\Q_p))\subset H^1(\Q_p,\pi_1(G))
$$
is a finite abelian group. Pick a set of representatives $R$, and replace $\|\cdot\|$ with
$$
\|v\|'=\max_{g \in R}\|g\cdot v\|.
$$
By construction, this modification $\|\cdot\|'$ is a $G(\Q_p)$-invariant norm on $\xi\otimes \pi$. $\square$

\subsection{The reductive case}

From the semisimple case, we derive the general reductive case.

\begin{df}
An irreducible admissible complex representation $\pi$ of $G(\Q_p)$ is essentially discrete series if a twist $\pi\otimes\nu$
is (unitary) discrete series, for some smooth character $\nu:G(\Q_p)\rightarrow \C^*$. The essentially discrete series representations over 
$\bar{\Q}_p$ are those of the form $\iota\pi$, for \underline{some} isomorphism $\iota:\C \rightarrow \bar{\Q}_p$. 
\end{df}

\noindent {\it{Remark}}. To put this definition (over $\bar{\Q}_p$) on more solid ground, we would like to know that we can in fact pick 
any $\iota$. In other words, whether any $\text{Aut}(\C)$-conjugate of an essentially discrete series representation is again essentially 
discrete series\footnote{Marko Tadic informs me that, at least for classical groups, this is known for generic representations.}. This is predicted by the local Langlands conjecture (the parameter does not map into a proper Levi). If $\sigma\in \Aut(\C)$, 
the matrix coefficients of $\sigma\pi$ are $\sigma$-conjugates of matrix coefficients of $\pi$. Hence, it is certainly true for supercuspidals, but
square integrability seems to be a problem. We should mention that at least it is known to be true for $\GL(n)$. Indeed the work of 
Bernstein-Zelevinsky shows that the essentially discrete series representations for $\GL(n)$ coincides with the generalized Steinberg 
representations $Q(\Delta)$, built from a segment $\Delta$ of supercuspidals, and $\sigma Q(\Delta)=Q(\sigma\Delta)$ in
a suitable (rational) normalization. See [Ku] for a nice exposition of the Langlands classification.

\begin{thm}
Let $G$ be a connected reductive group over $\Q_p$. Let $\xi$ be any irreducible algebraic representation of $G_{\bar{\Q}_p}$, and
$\pi$ be any essentially discrete series representation of $G(\Q_p)$ (both over $\bar{\Q}_p$). Then the locally 
algebraic representation $\xi\otimes \pi$ admits a $G(\Q_p)$-invariant norm if and only if its central character $\omega_{\xi}\cdot\omega_{\pi}$
is integral (that is, maps into $\bar{\Z}_p^{\times}$).
\end{thm}

\noindent {\it{Proof}}. The only if part is obvious. We assume $\omega_{\xi}\cdot\omega_{\pi}$
is integral, and seek a norm. The derived subgroup $G^{\text{der}}$ is semisimple, $Z_G \cap G^{\text{der}}$ is finite, and 
$$
1 \rightarrow Z_G \cap G^{\text{der}} \rightarrow Z_G\times G^{\text{der}}\rightarrow G \rightarrow 1.
$$
is exact. Here $Z_G$ is the full identity component of the center. The restriction $\xi|_{G^{\text{der}}}$ hence remains irreducible, and we will just 
write $\xi$. On the other hand, the restriction $\pi|_{G^{\text{der}}(\Q_p)}$ may not be, but it breaks up as a direct sum 
$$
\pi|_{G^{\text{der}}(\Q_p)}\simeq\tau_1\oplus \cdots \oplus \tau_r
$$
of discrete series representations $\tau_i$ of $G^{\text{der}}(\Q_p)$. For example, see [Ta, p. 381] and 
[Ta, p. 385]. By Lemma 5, there is a norm $\|\cdot\|_i$ on $\xi\otimes\tau_i$, invariant under $G^{\text{der}}(\Q_p)$. Their maximum defines a
$G^{\text{der}}(\Q_p)$-invariant norm $\|\cdot\|$ on $\xi\otimes \pi$, which is automatically 
$Z_G(\Q_p)$-invariant, by our assumption on the central character. 
$$
G(\Q_p)/Z_G(\Q_p)G^{\text{der}}(\Q_p)\subset H^1(\Q_p,Z_G \cap G^{\text{der}})
$$
is a finite abelian group. Pick representatives $R$, and replace $\|\cdot\|$ with
$$
\|v\|'=\max_{g \in R}\|g\cdot v\|.
$$
This is independent of $R$, and defines a $G(\Q_p)$-invariant norm on $\xi\otimes \pi$. $\square$

\medskip

\noindent Taking $G=\text{Res}_{F/\Q_p}\GL(n)$, for a finite extension $F/\Q_p$, yields:

\begin{cor}
Conjecture 5.5 in [BS] holds true. 
\end{cor}

\noindent {\it{Proof}}. As already mentioned, by Bernstein-Zelevinsky, the essentially discrete series representations of
$\GL(n)$ are precisely the generalized Steinberg representations. $\square$

%-----------------------------------------------------------------------------------------------------------

\noindent {\sc{Department of Mathematics, Princeton University, USA.}}

\noindent {\it{E-mail address}}: {\texttt{claus@princeton.edu}}

\end{document}